\documentclass{agtart_a}
\pdfoutput=1
\usepackage[all]{xy}

\title{Joins of DGA modules and sectional category}

\author[L Fern\'andez Su\'arez, P Ghienne, T Kahl and L
Vandembroucq]{Luc\'{\i}a Fern\'andez Su\'arez\newline Pierre
Ghienne\newline
Thomas Kahl\newline Lucile Vandembroucq}
\givenname{Luc\'{\i}a}
\surname{Fern\'andez Su\'arez}
\address{{\rm L Fern\'andez Su\'arez, T Kahl, L Vandembroucq:}\newline
Universidade do Minho\\
Centro de Matem\'atica\\
Campus de Gualtar\\\newline
4710-057 Braga\\
Portugal} 
\email{lfernandez@math.uminho.pt}

\givenname{Pierre}
\surname{Ghienne}
\address{{\rm P Ghienne:}\newline Universit\'e d'Artois\\
Facult\'e des Sciences Jean Perrin\\\newline
Laboratoire de Math\'ematiques\\ 
Rue Jean Souvraz - SP 18\\\newline
62307 Lens Cedex\\France} 
\email{ghienne@euler.univ-artois.fr}

\givenname{Thomas}
\surname{Kahl}
\email{kahl@math.uminho.pt}

\givenname{Lucile}
\surname{Vandembroucq}
\email{lucile@math.uminho.pt}



\received{1 November 2005}
\revised{11 January 2006}
\accepted{20 January 2006}
\proposed{}
\seconded{}
\published{24 February 2006}
\publishedonline{24 February 2006}

\volumenumber{6}
\issuenumber{}
\publicationyear{2006}
\papernumber{4}
\startpage{119}
\endpage{144}

\doi{}
\MR{}
\Zbl{}

\arxivreference{}
\arxivpassword{}

\AtBeginDocument{\let\bar\wbar}
\AtBeginDocument{\let\tilde\wtilde}

\makeatletter
\def\cnewtheorem#1[#2]#3{\newtheorem{#1}{#3}[section]
\expandafter\let\csname c@#1\endcsname\c@theor}

\makeautorefname{theor}{Theorem}

\newtheorem{theor}{Theorem}[section]
\cnewtheorem{lem}[theor]{Lemma}
\cnewtheorem{prop}[theor]{Proposition}
\cnewtheorem{cor}[theor]{Corollary}
\theoremstyle{definition}
\cnewtheorem{defin}[theor]{Definition}
\cnewtheorem{rem}[theor]{Remark}
\cnewtheorem{rems}[theor]{Remarks}

\cnewtheorem{exam}[theor]{Example}
\cnewtheorem{exs}[theor]{Examples}

\makeatother  

\newcommand{\fib}{\twoheadrightarrow }
\newcommand{\cof}{\rightarrowtail }
\newcommand{\we}{\smash{\rlap{\kern 6pt 
\raise 4pt\hbox{\footnotesize $\sim$}}}\longrightarrow}

\let\dem\proof
\let\findem\endproof

\newdir{ >}{{}*!/-7pt/\dir{>}}

\def\disfrac#1#2{{\displaystyle{\frac{#1}{#2}  }}}

\def\im{{\mbox{im}}}

\def\cat{{\rm{cat}\hskip1pt}}

\def\Mcat{{\rm{Mcat}\hskip1pt}}

\def\secat{{\rm{secat}\hskip1pt}}
\def\Msecat{{\rm{Msecat}\hskip1pt}}
\def\nil{{\rm{nil}\hskip1pt}}
\def\TC{{\rm{TC}\hskip1pt}}
\def\MTC{{\rm{MTC}\hskip1pt}}

\def\N{{\mathbb{N}}}

\begin{document}

\begin{asciiabstract}
We construct an explicit semifree model for the fiber join of two
fibrations p: E --> B and p': E' --> B from semifree models of p and
p'.  Using this model, we introduce a lower bound of the sectional
category of a fibration p which can be calculated from any Sullivan
model of p and which is closer to the sectional category of p than the
classical cohomological lower bound given by the nilpotency of the
kernel of p^*: H^*(B;\Q) --> H^*(E;\Q).  In the special case of the
evaluation fibration X^I --> X \times X we obtain a computable lower
bound of Farber's topological complexity TC(X).  We show that the
difference between this lower bound and the classical cohomological
lower bound can be arbitrarily large.
\end{asciiabstract}

\begin{htmlabstract}
We construct an explicit semifree model for the fiber join of two
fibrations p:E&#x2192;B and p':E'&#x2192;B from semifree
models of p and p'.  Using this model, we introduce a lower bound
of the sectional category of a fibration p which can be calculated
from any Sullivan model of p and which is closer to the sectional
category of p than the classical cohomological lower bound given by
the nilpotency of the kernel of p<sup>*</sup>:H<sup>*</sup>(B;<b>Q</b>)&#x2192;H<sup>*</sup>(E;<b>Q</b>).
In the special case of the evaluation fibration
X<sup>I</sup>&#x2192;X&times;X we obtain a computable lower bound
of Farber's topological complexity TC(X).  We show that the
difference between this lower bound and the classical cohomological
lower bound can be arbitrarily large.
\end{htmlabstract}

\begin{abstract}
We construct an explicit semifree model for the fiber join of two
fibrations $p\colon\, E \to B$ and $p'\colon\, E' \to B$ from semifree
models of $p$ and $p'$.  Using this model, we introduce a lower bound
of the sectional category of a fibration $p$ which can be calculated
from any Sullivan model of $p$ and which is closer to the sectional
category of $p$ than the classical cohomological lower bound given by
the nilpotency of the kernel of $p^*\colon\, H^*(B;\mathbb Q) \to
H^*(E;\mathbb Q)$.  In the special case of the evaluation fibration
$X^I \to X \times X$ we obtain a computable lower bound of Farber's
topological complexity $\mathrm{TC}(X)$.  We show that the difference
between this lower bound and the classical cohomological lower bound
can be arbitrarily large.
\end{abstract}

\maketitle

\section{Introduction}
The \textit{sectional category} of a fibration $p\co E\to B$, denoted
by $\secat p$, is the least integer $n$ such that the base space $B$
can be covered by $n+1$ open subspaces on each of which $p$ admits a
section. If no such $n$ exists one sets $\secat p = \infty$. This
homotopy invariant of a fibration has been introduced by A\,S Schwarz
\cite{Svarc} in the late 1950's as a generalization of the
Lusternik--Schnirelmann category of a space. The
Lusternik--Schnirelmann category of a space $X$, $\cat X$, is the
least integer $n$ such that $X$ can be covered by $n+1$ open subspaces
each of which is contractible in $X$ (if no such $n$ exists one sets
$\cat X = \infty$). If $X$ is a path-connected space with base point
$x_0$ and $PX$ is the space of paths beginning at $x_0$ then $\cat X$
is precisely the sectional category of the evaluation fibration
$ev_1\co PX \to X, \omega \mapsto \omega (1)$. References on
Lusternik--Schnirelmann category and sectional category are Schwarz
\cite{Svarc}, James \cite{James1,James2} and
Cornea--Lupton--Oprea--Tanr\'e \cite{CLOT}.

The concept of sectional category has been used to introduce measures
for the complexity of certain problems. S Smale \cite{Smale} (see also
\cite[sec. 9.4]{CLOT}) obtained results on the complexity of the
root-finding problem for algebraic equations in terms of sectional
category. Recently, M Farber \cite{Farber1,Farber2} defined the
\textit{topological complexity} of a space $X$, $\TC(X)$, to be the
sectional category of the evaluation fibration $ev_{0,1}\co X^I \to
X\times X, \omega \mapsto (\omega (0),\omega (1))$. This notion of
topological complexity plays an important role in the study of the
motion planning problem in robotics.

In spite of the simplicity of the definition, it is very hard to calculate the sectional category of a fibration $p \co  E \to B$ and therefore one will usually have to accept to work with approximations. For a surjective fibration one easily shows that $\secat p \leq \cat B$. Hence all upper bounds of $\cat B$, such as the dimension of $B$ or its cone-length, are upper bounds of $\secat p$ as well.
A classical cohomological lower bound of $\secat p$ is $\nil \ker p^*$, the nilpotency of the kernel of $p^* \co  H^*(B) \to H^*(E)$ (with respect to any coefficient ring), i.e.\ the least integer $n$ such that any ($n+1$)--fold cup product in $\ker p^*$ is trivial (cf.\ \cite[Section 9.3]{CLOT}). There are, of course, examples where $\nil \ker p^* = \secat p$ but in general the inequality $\nil \ker p^* \leq \secat p$ is strict.
As is showing the case of Lusternik--Schnirelmann category, that is, the special case where $p$ is the evaluation fibration $ev_1\co  PX \to X$,
the difference between the two numbers may actually be infinite.

A far better lower bound of $\secat p$ than $\nil \ker p^*$ (at least when the coefficient ring is $\Q$) is the rational sectional category
$\secat_0p$, i.e.\ the sectional category of a rationalization of $p$. In her thesis \cite{Agnese}, A Fass\`o Velenik gave a characterization
of $\secat_0p$ in terms of a Sullivan model of $p$. Unfortunately, concrete computations based on this characterization turn out to be rather difficult due to the complexity of the algebraic manipulations involved.
In the present article we introduce an approximation of $\secat p$ which is not as good as $\secat_0p$ in general but much easier to calculate. This approximation, which we denote by $\Msecat p$, is still a much better lower bound of $\secat p$ than $\nil \ker p^*$, if we consider coefficients in $\Q$.
Let us note here that we work over the field $\Q$ in the algebraic part of this article and that all spaces we consider are compactly generated Hausdorff spaces.

There is a classical equivalent definition of sectional category in terms of joins which is more appropriate for our purpose
than the original one. Denote by $\ast ^n_BE$ the $n$--fold fiber join of the fibration $p \co  E \to B$ and by $j^np \co  \ast^n_B E \to B$
the $n$th join map. If $B$ is normal then $\secat p \leq n$ if and only if $j^np$ has a section. We recall this fact and the join construction in section 2.
Let $A_{PL}$ denote Sullivan's functor of polynomial forms from spaces to commutative cochain algebras. Consider the morphism
$A_{PL}(j^np) \co  A_{PL}(B) \to A_{PL}(\ast^n_BE)$ as a morphism of $A_{PL}(B)$--modules. In section 5, we define the
invariant $\Msecat p$ to be the least integer $n$ for which $A_{PL}(j^np) = \phi \circ i$ where $\phi$ is a quasi-isomorphism of
$A_{PL}(B)$--modules and $i$ is a morphism of $A_{PL}(B)$--modules which admits a retraction of $A_{PL}(B)$--modules. We show
that $\nil \ker p^* \leq \Msecat p \leq \secat p$ for fibrations with a normal base space (cf.\ \fullref{ineq}). In the special case of the evaluation fibration $ev_1 \co  PX \to X$ of a simply connected space of finite type, $\Msecat$
coincides with the well-known invariant $\Mcat X$ (cf.\ \fullref{Mcat}) which in turn is known to be the rational category of $X$
(cf.\ Hess \cite{Hess}). The invariant $\Msecat$ generalizes the invariant $\Mcat$ hence in the same way as $\secat$ generalizes $\cat$.
The fact that $\Mcat$ is rational category does, however, not generalize to $\Msecat$ and $\Msecat$ does not in general equal rational
sectional category.

The computability of the invariant $\Msecat$ relies on an algebraic join construction which we develop in sections 3 and 4. Let $(A,d)$ be a commutative cochain algebra. In section 3, we define the join $(M,d) \ast_{(A,d)}(N,d)$ of two  $(A,d)$--semifree extensions $(M,d)$ and $(N,d)$ of $(A,d)$. This is an explicitly defined semifree extension of $(A,d)$. Moreover, if $(M,d)$ and $(N,d)$ are minimal semifree $(A,d)$--modules, so is $(M,d) \ast_{(A,d)}(N,d)$. Consider two fibrations $p\co  E \to B$ and $p'\co E' \to B$ between simply connected spaces of finite type and suppose that $\alpha \co  (A,d) \we A_{PL}(B)$ is a commutative cochain algebra model of the base space $B$ and that $(M,d)$ and $(N,d)$ are semifree extensions of $(A,d)$ such that there exist quasi-isomorphisms of $(A,d)$--modules $(M,d) \we A_{PL}(E)$ and $(N,d)\we A_{PL}(E')$ which extend $A_{PL}(p) \circ \alpha$ and $A_{PL}(p') \circ \alpha$. We establish in section 4 that the inclusion $(A,d) \to (M,d) \ast_{(A,d)}(N,d)$ is an $(A,d)$--module model of the topological join map $E\ast_BE'\to B$. Iterating the join construction, we define the $n$--fold join $\ast^n_{(A,d)}(M,d)$ of $(M,d)$ and obtain an explicit $(A,d)$--module model of the $n$th join map $j^np \co  \ast^n_BE \to B$. The number $\Msecat p$ is then the least $n$ such that the inclusion $(A,d) \to \ast^n_{(A,d)}(M,d)$ admits a retraction of $(A,d)$--modules (cf.\ \fullref{Mjoin}). Through this result one obtains an effective method to compute the invariant $\Msecat p$ from a Sullivan model of $p$.

As an example we consider Farber's topological complexity $\TC$. Let $X$ be a simply connected space of finite type with
Sullivan model $(\Lambda V,d)$. There is a well-known explicit minimal model of the evaluation fibration $ev_{0,1} \co  X^I \to X\times X$
which can be determined from $(\Lambda V,d)$. This model and the algebraic join construction permit one to calculate the invariant
$\MTC(X) = \Msecat ev_{0,1}$ which is a lower bound of $\TC(X)$. Note that since $ev_{0,1}$ is the mapping path fibration associated to the diagonal map
$X \to X\times X$, $ev_{0,1}^*$ can be identified with the cup product $\cup \co  H^*(X)\otimes H^*(X) \to H^*(X)$. If $X$ is a formal space,
i.e.\ a space whose rational homotopy type is entirely determined by its cohomology algebra, one has
$\MTC(X) = \nil \ker \cup$. But already for the simplest example of a non-formal space, one calculates that $\MTC(X) =3$ and
$\nil \ker \cup = 2$. We show finally that the difference between the two lower bounds is unbounded.

\section{Sectional category and joins}
Recall from the introduction that the category of spaces in which we shall work throughout this article is the category of compactly generated Hausdorff spaces. All categorical constructions (products, pullbacks etc.) are carried out in this category.

In this section we recall the link between joins and the sectional category mentioned in the introduction.

\begin{defin}\rm 
The \textit{(fiber) join} of two maps $p \co  E \to B$ and $p'\co E'\to B$, denoted by $E\ast _BE'$, is the double mapping cylinder of the projections $E\times _BE' \to E$ and $E\times _BE' \to E'$, i.e.\ the quotient space $((E\times _BE')\times I \amalg E \amalg E')/\sim$ where $(e,e',0) \sim e$, $(e,e',1) \sim e'$. The \textit{join map} of $p$ and $p'$ is the map $j_{p,p'} \co  E\ast_BE' \to B$ defined by $j_{p,p'}([e,e',t]) = p(e) = p'(e')$, $j_{p,p'}([e]) = p(e)$, and $j_{p,p'}([e']) = p'(e')$. The \textit{$n$--fold join} and the \textit{$n$th join map} of $p$ are iteratively defined by $\ast^0_BE = E $, $\ast^n_BE = (\ast^{n-1}_BE) \ast_BE$, $j^0p = p$, and $j^np = j_{j^{n-1}p,p}$.
\end{defin}

\begin{theor}\label{joinsecat}
Let $p\co E\to B$ be a fibration. If $B$ is normal then $\secat p \leq n$ if and only if $j^np$ has a section.
\end{theor}

\dem The result is well-known, at least when $B$ is paracompact (cf.\
James \cite{James1}). We include a short proof for the convenience of
the reader.

Suppose first that $\secat p \leq n$. We show by induction that for each $0 \leq m\leq n$ there exists an open cover $U_0, \dots, U_{n-m}$ of $B$ such that $j^{m}p$ has a section on $U_0$ and $p$ has a section on each of the remaining $U_i$. For $m= 0$ this is just the hypothesis that $\secat p \leq n$. Suppose that the assertion holds for $0 \leq m < n$. Then there exists an open cover $U_0, \dots, U_{n-m}$ of $B$, a section $\sigma _0 \co  U_0 \to \ast^{m}_BE$ of $j^mp$, and sections $\sigma _i \co  U_i \to E$ of $p$ $(1 \leq i \leq n-m)$.  Since $B$ is normal, there exist open covers $V_0,\dots ,V_{n-m}$ and $W_0,\dots , W_{n-m}$ of $B$ such that $\bar V_i \subset W_i \subset \bar W_i \subset U_i$. Set $A_0 = \bar V_0 \cap (B\backslash W_1)$, $A_1 = \bar V_1 \cap (B\backslash W_0)$, and $A_2 = \bar W_0 \cap \bar W_1 \cap (\bar V_0 \cup \bar V_1)$. Then $A_0$, $A_1$, and $A_2$ are closed subspaces of $B$, $A_0\cup A_1 \cup A_2 = \bar V_0 \cup \bar V_1$, and $A_0 \cap A_1 = \emptyset$. Since $B$ is normal, by Urysohn's Lemma, there exists a continuous map $\phi \co  B \to I$ such that $\phi (A_0) \subset \{ 0\}$ and $\phi (A_1) \subset \{ 1\}$. Define a section $\sigma $ of $j^{m+1}p$ on $\bar V_0 \cup \bar V_1$ by 
$$\sigma (x) = \left\{ \begin{array}{lll}
\left[\sigma _0(x)\right], &  x\in A_0,\\
\left[\sigma_1(x)\right], &   x\in A_1,\\
\left[\sigma_0(x), \sigma_1(x),\phi(x)\right],  &   x\in A_2.
\end{array}\right.$$
Consider the open cover $O_0,\dots O_{n-m-1}$ of $B$ given by $O_0 = V_0\cup V_1$ and $O_i = U_{i+1}$, $i = 1,\dots ,n-m-1$. 
On $O_0$, $\sigma$ is a section of $j^{m+1}p$. On each of the remaining $O_i$, $p$ has a section by hypothesis. 
This terminates the inductive step.

Suppose now that $j^np$ has a section $s \co  B \to \ast^n_BE$. By \fullref{secatq} below, $\ast^n_BE$ can be covered by $n+1$ open subspaces $U_0, \dots, U_n$ on each of which the projection $\bar p_n \co  (\ast^n_BE) \times _BE \to \ast^n_BE$ has a section. The inverse images $s^{-1}(U_i)$ form a cover of $B$ by open subspaces on each of which $p$ has a section. Therefore $\secat p \leq n$.
\findem

\begin{rem}\rm 
If $p\co  E \to B$ and $p'\co E'\to B$ are fibrations, so is the join map $j_{p,p'} \co  E\ast_BE' \to B$. Indeed, if $\lambda \co  E\times _BB^I \to E^I$ and $\lambda' \co  E'\times _BB^I \to E'^I$ are lifting maps for $p$ and $p'$ then 
a lifting map $\phi \co  (E\ast _BE')\times _BB^I \to (E\ast_B E')^I$ for $j_{p,p'}$ is given by
$\phi ([e,e',t],\omega) (s) = [\lambda (e,\omega )(s),\lambda'(e',\omega )(s),t]$, $\phi ([e],\omega) (s) = [\lambda (e,\omega )(s)]$, and $\phi ([e'],\omega) (s) = [\lambda' (e',\omega )(s)]$. Note that $\phi$ is continuous since we are working with compactly generated spaces. It follows, by induction, that the $n$th join map of a fibration is again a fibration and hence that it has a section if and only if it has a homotopy section.
\end{rem}

In the proof of \fullref{joinsecat} we used the following lemma. We shall need this lemma again in the proof of the inequality $\Msecat p \geq \nil \ker p^*$ (cf.\ \fullref{ineq}).

\begin{lem}\label{secatq}
Consider a fibration $p\co  E \to B$ and form the pullback diagram:
$$\disablesubscriptcorrection\xymatrix{
(\ast^n_BE)\times_BE \ar[r] \ar[d]_{\bar p_n} 
& E\ar [d]^{p}
\\ 
\ast^n_BE \ar[r]_{j^np}
& B
}$$
Then $\secat \bar p_n \leq n$. 
\end{lem}

\dem
We proceed by induction. For $n = 0$, the map $E \to E\times _BE$, $e \mapsto (e,e)$ is a section of $\bar p_n$. Suppose that $n>0$ and that the assertion holds for $n-1$. The spaces $E$ and $\ast^{n-1}_BE$ are embedded as closed subspaces in $\ast^n_BE$ and there are canonical projections $\pi \co  \ast^n_BE \backslash \ast^{n-1}_BE \to E$ and $\tilde \pi \co  \ast^n_BE \backslash E \to \ast^{n-1}_BE$. Let $U_0$ be the open subspace $\ast^n_BE \backslash \ast^{n-1}_BE$ of $\ast^n_BE$. We have $j^np|_{U_0} = p\pi$. The inductive hypothesis implies that $\ast^{n-1}_BE$ can be covered by $n$ open subspaces $V_1, \dots , V_{n}$ such that each restriction of the join map $j^{n-1}p|_{V_i}\co  V_i \to B$ factors through $p$. Consider the open subspaces $U_i = \tilde \pi^{-1}(V_i)$ of $\ast^n_BE$. The $n+1$ open subspaces $U_0, U_1 \dots, U_n$ of $\ast^n_BE$ cover $\ast^n_BE$. The restriction of the join map $j^np$ to any of these open subspaces factors through $p$. Therefore the projection $\bar p_n \co  (\ast^n_BE)\times _BE \to \ast^n_BE$ has a section on each $U_i$. This shows that $\secat \bar p_n \leq n$.
\findem

\section{Joins of semifree modules}

The purpose of this section is to define joins of semifree extensions of a commutative cochain algebra. Recall that we are working over $\Q$. All graded vector spaces we consider will be $\Z$--graded with upper degree and all differential vector spaces will be cochain complexes, i.e.\ the differential raises the upper degree by one. The \textit{$n$th suspension} $s^{-n}V$ of a graded vector space $V$ is defined by $(s^{-n}V)^i = V^{i-n}$. 

\begin{defin}\rm
Let $(A,d)$ be a differential algebra. A \textit{semifree extension} of an $(A,d)$--module $(M,d)$ is an $(A,d)$--module of the form $(M\oplus A\otimes X,d)$ where the action is the one of the direct sum, the differential on $M$ is the differential of $(M,d)$, and $X$ admits a direct sum decomposition $X = \bigoplus \limits_{i=0}^{\infty} X_i$ such that $d(X_0) \subset M$ and $d(X_n) \subset M \oplus A\otimes (\bigoplus \limits_{i=0}^{n-1} X_{i})$ for $n \geq 1$. A \textit{semifree $(A,d)$--module} is a semifree extension of the trivial $(A,d)$--module $0$. 
\end{defin}

For the remainder of this section we fix a commutative cochain algebra $(A,d)$ and two semifree extensions $(M,d) = (A \oplus A\otimes X,d)$ and $(N,d) =  (A \oplus A\otimes Y,d)$ of $(A,d)$. We define the join $(M,d) \ast _{(A,d)}(N,d)$ of $(M,d)$ and $(N,d)$ which will again be a semifree extension of $(A,d)$. Forgetting the differential, $(M,d)\ast_{(A,d)}(N,d)$ is the free graded $A$--module $A \oplus A\otimes s^{-1}X\otimes Y$.
In order to define the differential of $(M,d)\ast_{(A,d)}(N,d)$, we decompose the differential in $(M,d)$ of an element $m \in M$ as  
$$dm = d_0m + d_+m$$ where 
$d_0m \in A$ and $d_+m \in A\otimes X$. Using the same notation, we decompose the differential in $(N,d)$ of an element $n \in N$. Consider elements $x\in X$ and $y \in Y$ and write
$$d_+x = \sum_i a_i\otimes x_i, \quad \mbox{and} \quad d_+y = \sum_j b_j\otimes y_j.$$
The differential of the element $s^{-1}x\otimes y$ in $(M,d)\ast_{(A,d)}(N,d)$ is then defined by
\begin{eqnarray*}
d(s^{-1}x\otimes y) &=& (-1)^{|x|}d_0xd_0y + \sum_i (-1)^{|a_i|+1}a_i\otimes s^{-1}x_i\otimes y\\
& &+ \sum_j (-1)^{(|x|+1)(|b_j|+1)}b_j\otimes s^{-1}x\otimes y_j.
\end{eqnarray*}
We extend this differential to the whole join $(M,d)\ast_{(A,d)}(N,d)$ by setting
$$d(a\otimes s^{-1}x\otimes y) = da\otimes s^{-1}x\otimes y + (-1)^{|a|}a\cdot d(s^{-1}x\otimes y).$$
\fullref{d^2=0} below assures that $d$ is indeed an $(A,d)$--module differential in\break $(M,d)\ast_{(A,d)}(N,d)$. 
It is an easy exercise to check that $(M,d)\ast_{(A,d)}(N,d)$ is a semifree extension of $(A,d)$. Moreover, if $(A,d)$ is augmented and $(M,d)$ and $(N,d)$ are minimal semifree $(A,d)$--modules, i.e.\ the differentials in  $\Q\otimes_{(A,d)}(M,d)$ and $\Q\otimes_{(A,d)}(N,d)$ are zero, then  $(M,d)\ast_{(A,d)}(N,d)$ is also minimal. 

\begin{prop} \label{d^2=0}
$d^2(s^{-1}x\otimes y)=0$.
\end{prop}

\dem 
Write $d_+x_i = \sum_k a_{ik}\otimes x_{ik}$ and $d_+y_j = \sum_l b_{jl}\otimes y_{jl}$. Since 
\begin{eqnarray*}
0& =& d^2x\\
&=& d(d_0x + \sum_i a_i\otimes x_i)\\
&=& dd_0x + \sum_i da_i\otimes x_i + \sum _i (-1)^{|a_i|}a_id_0x_i + \sum _i (-1)^{|a_i|}a_id_+x_i,
\end{eqnarray*}
we have $dd_0x = - \sum _i (-1)^{|a_i|}a_id_0x_i = \sum _i (-1)^{|a_i|+1}a_id_0x_i$ and 
$$\sum_i da_i\otimes x_i = -\sum _i (-1)^{|a_i|}a_id_+x_i = \sum_{i,k} (-1)^{|a_i|+1}a_ia_{ik}\otimes x_{ik}.$$
Similarly, 
$dd_0y = \sum _j (-1)^{|b_j|+1}b_jd_0y_j$ and 
$$\sum_j db_j\otimes y_j = \sum_{j,l} (-1)^{|b_j|+1}b_jb_{jl}\otimes y_{jl}.$$
Use $\sigma$ to denote the isomorphism 
$$A\otimes X\otimes Y \to A\otimes s^{-1}X\otimes Y,\quad a\otimes x\otimes y \mapsto (-1)^{|a|}a\otimes s^{-1}x\otimes y$$
and $T$ to denote the isomorphism
$$A\otimes X \to X\otimes A,\quad a\otimes x \mapsto (-1)^{|a||x|}x\otimes a.$$
Applying $\sigma$ to the identity 
$\sum_i da_i \otimes x_i\otimes y = \sum_{i,k} (-1)^{|a_{i}|+1}a_ia_{ik}\otimes x_{ik}\otimes y$,
one obtains the identity
$$\sum_i (-1)^{|a_i|+1}da_i \otimes s^{-1}x_i\otimes y = \sum_{i,k} (-1)^{|a_{ik}|+1}a_ia_{ik}\otimes s^{-1}x_{ik}\otimes y.$$
Applying $\sigma \circ (T\otimes id_Y)$ to the identity
$$\sum_j x\otimes db_j \otimes y_j = \sum_{j,l} (-1)^{|b_j| +1}x\otimes b_jb_{jl}\otimes  y_{jl},$$
one obtains the identity
\begin{align*}
\sum_j (-1)^{|x| + |b_j| + |x||b_j|}&db_j \otimes s^{-1}x\otimes y_j\\
=&\sum_{j,l} (-1)^{|x||b_j|+|x||b_{jl}|+ |b_{jl}|}b_jb_{jl}\otimes s^{-1}x\otimes y_{jl}.
\end{align*}
Using the different formulae above one easily verifies that $d^2(s^{-1}x\otimes y) =0$.
\findem

\begin{defin}\rm
The \textit{$n$--fold join} of  $(M,d)$ is iteratively defined by  $\ast_{(A,d)}^0(M,d) = (M,d)$ and  $\ast^n_{(A,d)}(M,d) = (\ast^{n-1}_{(A,d)}(M,d))\ast_{(A,d)}(M,d)$. 
\end{defin}

\begin{rems} \rm \label{iteratedjoin}
(i)\qua Note that $\ast^n_{(A,d)}(M,d)$ is a semifree extension of $(A,d)$. Moreover, if $(M,d)$ is a minimal semifree $(A,d)$--module then $\ast^n_{(A,d)}(M,d)$ is a minimal semifree $(A,d)$--module as well. 

(ii)\qua We have $$\ast^n_{(A,d)}(M,d) = (A \oplus A\otimes s^{-n}X^{\otimes n+1},d).$$ 
Consider elements $x_0, \dots ,x_n \in X$ and write $d_+x_i = \sum \limits_{j_i}a_{ij_i}\otimes x_{ij_i}$. An easy induction shows that
\begin{eqnarray*}
\lefteqn{d(s^{-n}x_0\otimes \cdots \otimes x_n) = (-1)^{\sum \limits_{k=1}^{n}(k|x_{n-k}| + k -1)} d_0x_0\cdot \cdots \cdot d_0x_n}\\
&+ & \sum \limits _{i=0}^{n} \sum \limits_{j_i}(-1)^{(|a_{ij_i}| +1)(|x_0| + \cdots +|x_{i-1}| + n)}a_{ij_i}\otimes s^{-n}x_0 \otimes \cdots \otimes x_{ij_i}\otimes  \cdots \otimes x_n.
\end{eqnarray*}

(iii)\qua Consider a morphism of commutative cochain algebras $(A,d) \to (B,d)$ and the $(B,d)$--semifree extension of $(B,d)$ defined by $(N,d) = (B,d)\otimes _{(A,d)}(M,d)$. The formula for the differential given in (ii) shows that  $$\ast^n_{(B,d)}(N,d) = (B,d)\otimes _{(A,d)}\ast^n_{(A,d)}(M,d).$$
\end{rems}

\section{Topological versus algebraic joins}
Our goal in this section is to show that the algebraic joins of the preceding section model topological joins.
In this and the following sections we make frequent use of the homotopy theory of modules over a DGA and, in particular, of the
following well-known result:

\begin{theor}\label{Model}
Let $(A,d)$ be a differential algebra. The category of $(A,d)$--modules is a proper closed model category where weak
equivalences are  quasi-isomorphisms, fibrations are surjective morphisms, and cofibrations are morphims having the left
lifting property with respect to surjective  quasi-iso\-morphisms. A morphism is a cofibration if and only if it is a retract of the inclusion of a semifree extension.
\end{theor}

We refer the reader to Goerss and Jardine \cite{GoerssJardine} for the axioms of closed model categories. A closed model category is called \textit{proper} if
the class of weak equivalences is closed under base change along fibrations and cobase change along cofibrations. As is customary
we denote weak equivalences by $\we$, fibrations by $\fib$, and cofibrations by $\cof$. For the convenience of the reader we include
the following proof.

\medskip\textbf{Proof of \fullref{Model}}\qua We do not use the fact that we are working over $\Q$ and the proof works for an arbitrary commutative ground ring. We first show that
inclusions of semifree extensions are cofibrations. Consider a semifree extension
$(M\oplus A\otimes (\bigoplus \limits_{i=0}^{\infty} X_i),d)$ of an $(A,d)$--module $(M,d)$ and a commutative diagram of $(A,d)$--modules
$$\disablesubscriptcorrection\xymatrix{
(M,d) \ar[d]_{i} \ar[r]^{f}
& (P,d) \ar@{->>} [d]^{p}_{\hbox{\small$\sim$}}
\\
(M\oplus A\otimes (\bigoplus \limits_{i=0}^{\infty} X_i),d) \ar[r]_-{g}
& (Q,d)
}$$
where $i$ is the inclusion. Suppose inductively that we have constructed a lifting $\lambda$ for the diagram up to
$(M\oplus A\otimes (\bigoplus \limits_{i=0}^{n} X_i),d)$. Let $\mathcal{ B} \subset X_{n+1}$ be a basis and $x \in \mathcal{B}$. Then
$\lambda (dx)$ is defined and $d\lambda (dx) = 0$. Since $p$ is surjective, there exists an element $\xi \in P$ such that $p(\xi )= g(x)$.
Then $\lambda (dx) - d\xi$ is a cocycle in $\ker p$. Since $p$ is a surjective quasi-isomorphism, $\ker p$ is acyclic and there exists an
element $y \in \ker p$ such that $dy = \lambda (dx) - d\xi$. Set $\lambda (x) = \xi + y$. This defines $\lambda$ in  $(M\oplus A\otimes (\bigoplus \limits_{i=0}^{n+1} X_i),d)$. It follows that a lifting exists and hence that $i$ is a cofibration. Axiom CM1 (existence of finite limits and colimits) follows from the fact that the category of modules over the ground ring is complete and cocomplete. The fact that the quasi-isomorphisms have the ``$2=3$" property (CM2) and are closed under retracts follows from the corresponding properties of isomorphisms. The fibrations are closed under retracts because surjective maps are closed under retracts. It is a general fact that any class of morphisms in a category which is defined by having the left lifting property with respect to another class of morphisms is closed under retracts. Therefore the class of cofibrations is closed under retracts and CM3 holds. We check the factorization axiom CM5. Consider a morphism $f\co  (M,d) \to (N,d)$ of $(A,d)$--modules.
Consider the acyclic semifree $(A,d)$--module $(A\otimes (\tilde N \oplus s^{-1}\tilde N), \delta)$ where $\tilde N = \bigoplus \limits_{n\in N} \Q\cdot n$, $\delta n = s^{-1}n$ and
$\delta s^{-1}n = 0$. Then the inclusion $i\co  (M,d) \to  (M,d) \oplus (A\otimes (\tilde N\oplus s^{-1}\tilde N), \delta)$ is both a cofibration
and a quasi-isomorphism. Let $p \co  (M,d) \oplus (A\otimes (\tilde N\oplus s^{-1}\tilde N), \delta) \to (N,d)$ be the morphism of $(A,d)$-modules defined by
$p(m) = f(m)$, $p(n) = n$, and $p(s^{-1}n) = dn$. Obviously, $p$ is surjective and $f = p\circ i$. This shows one part of CM5.
In the proof of F\'elix--Halperin--Thomas, \cite[2.1(i)]{FHTDGA}, it is shown that there is a factorization $f = p\circ i$
where $p$ is a surjective quasi-isomorphism and $i$ is the inclusion of a semifree extension. This shows the other part of CM5.
We verify the lifting axiom CM4. One of the lifting properties is the definition of cofibrations. For the other one consider a
commutative diagram of $(A,d)$--modules:
$$\disablesubscriptcorrection\xymatrix{
(M,d) \ar@{ >->}[d]^{i}_{\hbox{\small$\sim$}} \ar[r]^{f}
& (P,d) \ar@{->>} [d]^{p}
\\
(N,d) \ar[r]_{g}
& (Q,d)
}$$
Form the pullback $(A,d)$--module $(N,d)\times _{(Q,d)}(P,d)$. Since $p$ is surjective, so is its base extension
$\bar p \co  (N,d)\times _{(Q,d)}(P,d) \to (N,d)$. Choose a factorization of the canonical morphism $(i,f)\co  (M,d) \to (N,d)\times _{(Q,d)}(P,d)$
in a quasi-isomorphism $j\co  (M,d)\we (R,d)$ and a surjective morphism $r\co  (R,d) \to (N,d)\times _{(Q,d)}(P,d)$. The
composite $\bar p\circ r$ is a surjective quasi-isomorphism. Since $i$ is a cofibration, there exists a section $s$ of $\bar p \circ r$ such
that $s\circ i = j$. Let $\bar g \co  (N,d)\times _{(Q,d)}(P,d) \to (P,d)$  be the base extension of $g$. Then the composite $\bar g \circ r \circ s$ is a
lifting for the above square. It follows that the category of $(A,d)$--modules is a closed model category.

We have seen that an inclusion of a semifree extension is a cofibration. By CM3, any retract of an inclusion of a semifree extension is a cofibration.  Let $i$ be a cofibration and $i = p \circ j$ be a factorization such that $j$ is the inclusion of a semifree extension and $p$ is a surjective quasi-isomorphism. We have already mentioned that such a factorization exists. By CM4, there exists a section $s$ of $p$ such that $s\circ i =j$. This implies that $i$ is a retract of $j$. Thus a morphism is a cofibration if and only if it is a retract of the inclusion of a semifree extension. In particular, any cofibration is injective.
Therefore we may use the $5$--lemma to show that the cobase extension of a weak equivalence along a cofibration is a weak equivalence.
Since, by definition, fibrations are surjetive, the $5$--lemma implies that base extension of a weak equivalence along a fibration is a weak equivalence.
It follows that the closed model category of $(A,d)$--modules is proper.
\findem

Consider two fibrations $p\co  E \to B$ and $p'\co  E' \to B$ where $B$, $E$, and $E'$ are simply connected spaces of finite type. Simply connected spaces are understood to be non-empty. A space is said to be of finite type if it has finite dimensional rational homology in every dimension. Let $A_{PL}$ be Sullivan's functor from spaces to commutative cochain algebras. Fix a commutative cochain algebra model $\alpha \co  (A,d) \we A_{PL}(B)$. For any continuous map $f\co  S \to B$ the morphism of commutative cochain algebras $A_{PL}(f) \circ \alpha$ induces an
$(A,d)$--module structure on $A_{PL}(S)$ such that $A_{PL}(f) \circ \alpha$ is a morphism of $(A,d)$--modules. Let $(M,d)= (A \oplus A\otimes X,d)$ and $(N,d)= (A \oplus A\otimes Y,d)$ be semifree extensions of $(A,d)$ such that there exist quasi-isomorphisms of $(A,d)$--modules $(M,d) \we A_{PL}(E)$ and $(N,d)\we A_{PL}(E')$ which extend $A_{PL}(p) \circ \alpha$ and $A_{PL}(p') \circ \alpha$. As in the preceding section we write $d = d_0 + d_+$ for the differentials of $(M,d)$ and $(N,d)$.

\begin{theor}\label{joinmodel}
{\rm(i)}\qua The morphism of $(A,d)$--modules $A_{PL}(j_{p,p'}) \circ \alpha$ extends to a quasi-isomorphism of $(A,d)$--modules
$(M,d)\ast_{(A,d)}(N,d) \we A_{PL}(E\ast_BE')$.

{\rm(ii)}\qua The morphism of $(A,d)$--modules $A_{PL}(j^np)\circ \alpha$ extends to a quasi-isomor\-phism of $(A,d)$--modules $\ast_{(A,d)}^n(M,d) \we  A_{PL}(\ast _B^nE)$.
\end{theor}

\dem
The second part follows from the first by induction. The proof of (i) is divided in $3$ steps.
\subsection*{Step 1:\qua A model of the pullback}
Choose Sullivan models $\psi \co  (A\otimes \Lambda V,d) \we A_{PL}(E)$ and $\psi ' \co  (A\otimes \Lambda V',d)  \we A_{PL}(E')$ of $A_{PL}(p) \circ \alpha$ and $A_{PL}(p') \circ \alpha$. Since the inclusions $(A,d) \to (M,d)$ and $(A,d) \to (N,d)$ are cofibrations, by the lifting lemma (Baues \cite[II.1.11]{Baues}), there exist quasi-isomorphisms of $(A,d)$--modules $h \co  (M,d)\we (A\otimes \Lambda V,d)$ and $h'\co (N,d) \we (A\otimes \Lambda V',d)$ which extend the inclusions of $(A,d)$. Form the pushout of commutative cochain algebras:
$$\disablesubscriptcorrection\xymatrix{
(A,d) \ar[r]\ar[d] 
& (A\otimes \Lambda V,d) \ar [d] 
\\ 
(A\otimes \Lambda V',d)\ar[r] 
& (A\otimes \Lambda V,d)\otimes_{(A,d)}(A\otimes \Lambda V',d)
}$$
It follows from F\'elix--Halperin--Thomas \cite[15(c)]{RHT} that the morphisms $A_{PL}(pr_E) \circ \psi$ and $A_{PL}(pr_E') \circ \psi'$, where $pr_E \co  E \times _BE' \to E$ and $pr_{E'}\co E\times _BE' \to E'$ are the projections, induce a quasi-isomorphism of commutative cochain algebras $$(A\otimes \Lambda V,d)\otimes_{(A,d)}(A\otimes \Lambda V',d)\to A_{PL}(E\times _BE').$$  
By \cite[6.7]{RHT}, since $h$ and $h'$ are quasi-isomorphisms between semifree $(A,d)$--modules, the morphism 
$$h\otimes _Ah'\co  (M,d)\otimes _{(A,d)}(N,d) \to (A\otimes \Lambda V,d)\otimes_{(A,d)}(A\otimes \Lambda V',d)$$
is a quasi-isomorphism. Since $A$ is commutative, $(M,d)\otimes _{(A,d)}(N,d)$ is an $(A,d)$--module and $h\otimes _Ah'$ a quasi-isomorphism of $(A,d)$--modules. Note that 
$$(M,d)\otimes _{(A,d)}(N,d) = (M\otimes _AN,d) = (A \oplus A\otimes (X \oplus Y \oplus X\otimes Y),d)$$ 
contains both $(M,d)$ and $(N,d)$ as sub $(A,d)$--modules. Note also that if $d_+x =  \sum_i a_i\otimes x_i$  and $d_+y =  \sum _jb_j\otimes y_j$ then the differential of $x\otimes y$ in  $(M,d)\otimes _{(A,d)}(N,d)$ is given by 
\begin{eqnarray*}
d(x\otimes y) &=& d_0x\otimes y + \sum_i a_i\otimes x_i \otimes y + (-1)^{|x||y|}d_0y\otimes x\\
&& + \sum_j (-1)^{|x|(|b_j|+1)}b_j\otimes x \otimes y_j.
\end{eqnarray*}
We have obtained the following commutative diagram of $(A,d)$--modules:
$$\disablesubscriptcorrection
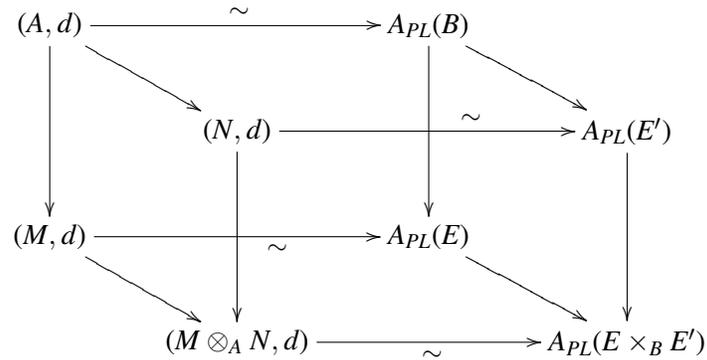
\xymatrix{
(A,d) \ar[rr]^{\hbox{\small$\sim$}} \ar[dr] \ar[dd]  &
& A_{PL}(B)\ar[rd] \ar[dd] &  
\\ & 
(N,d) \ar[rr]^(.6){\hbox{\small$\sim$}} \ar[dd] & 
& A_{PL}(E')\ar[dd]
\\ 
(M,d) \ar[rr]_(.6){\hbox{\small$\sim$}} \ar[dr] & 
& A_{PL}(E)\ar[rd] & 
\\ & 
(M\otimes _{A}N,d) \ar[rr]_{\hbox{\small$\sim$}} &
& A_{PL}(E\times_BE')
}$$

\subsection*{Step 2:\qua A model of the join map}
Consider the mapping cylinder factorization of the projection $pr_E' \co  E\times _BE' \to E'$ in a cofibration $\iota \co  E\times _BE' \to Z$ and a homotopy equivalence $\rho\co  Z \to E'$. We have the following pushout:
$$\disablesubscriptcorrection\xymatrix{
E\times _BE'\ar[d]_{pr_E} \ar[r]^{\iota} 
& Z\ar [d] 
\\ 
E\ar[r] 
& E\ast_BE'
}$$
\eject
Let $\nu$ be the inclusion $(N,d) \to (M\otimes _AN,d)$. We construct the mapping path factorization of $\nu$. Consider the $(A,d)$--module 
$$(Q,D) = (M \otimes_{A}N \oplus  N \oplus s^{-1}M\otimes_{A}N,D)$$
where the action on $M \otimes_{A}N \oplus  N$ is the one of the direct sum, $a\cdot s^{-1}w =$\break $(-1)^{|a|}s^{-1}aw$, 
and the differential is given by
\begin{eqnarray*} 
D(m\otimes _An)  &=& d(m\otimes _An) + s^{-1}m\otimes _An,\\ 
Dn &=& dn + s^{-1}\nu(n),\\  
Ds^{-1}w &=& -s^{-1}dw.
\end{eqnarray*}
Let $i \co  (N,d) \to (Q,D)$ be the injection defined by $i(n) = \nu(n) - n$. One easily checks that this is a morphism of $(A,d)$--modules. We show that $i$ is both a cofibration and a quasi-isomorphism. Set $U = \Q \oplus X \oplus Y \oplus X\otimes Y$ where the elements of $\Q$ have degree $0$. Then $M\otimes _AN$ is the free graded $A$--module $A\otimes U$. Consider the acyclic semifree $(A,d)$--module $(A\otimes (U \oplus s^{-1}U), \delta)$ where $\delta u = s^{-1}u$ and $\delta s^{-1}u = 0$. Then the inclusion 
$$(N,d) \to  (N,d) \oplus (A\otimes (U\oplus s^{-1}U), \delta)$$ is both a cofibration and a quasi-isomorphism.   
Consider the isomorphism of $(A,d)$--modules $\Phi\co  (N,d) \oplus (A\otimes (U \oplus s^{-1}U), \delta) \to (Q,D)$ defined by $\Phi(n) = \nu (n) - n$, $\Phi (u) = u$, and $\Phi(s^{-1}u)= du + s^{-1}u$. Since $i$ is the restriction of $\Phi$ to $(N,d)$, it is both a cofibration and a quasi-isomorphism. Let $\pi \co  (Q,D) \to (M\otimes _{A}N,d)$ be the obvious projection. Then $\pi$ is a surjective morphism of $(A,d)$--modules and $\pi \circ i = \nu$.

Form the following commutative diagram of $(A,d)$--modules:
$$\disablesubscriptcorrection\xymatrix{
(N,d) \ar[r]_{\hbox{\small$\sim$}} \ar@{ >->}[d]^{i}_{\hbox{\small$\sim$}}
& A_{PL}(E') \ar[r]^{A_{PL}(\rho)}_{\hbox{\small$\sim$}}
& A_{PL}(Z) \ar [d]^{A_{PL}(\iota)} 
\\ 
(Q,D)\ar[r]_-{\pi}
& (M\otimes_AN,d) \ar[r]^-{\hbox{\small$\sim$}}
& A_{PL}(E\times_BE')
}$$
Since $A_{PL}(\iota )$ is surjective, there exists a lifting $\lambda \co  (Q,D)\to A_{PL}(Z)$ making the diagram commutative. Note that $\lambda $ is automatically a quasi-isomorphism. Consider the following commutative cube of $(A,d)$-modules: 
$$\disablesubscriptcorrection
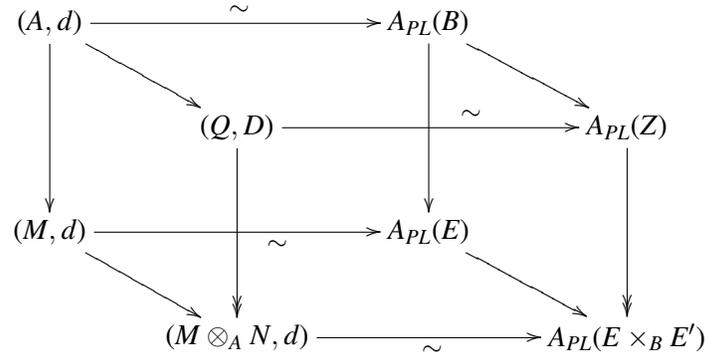
\xymatrix{
(A,d) \ar[rr]^{\hbox{\small$\sim$}} \ar[dr] \ar[dd]  &
& A_{PL}(B)\ar[rd] \ar[dd] &  
\\ & 
(Q,D) \ar[rr]^(.6){\hbox{\small$\sim$}} \ar@{->>}[dd] & 
& A_{PL}(Z)\ar@{->>}[dd]
\\ 
(M,d) \ar[rr]_(.6){\hbox{\small$\sim$}} \ar[dr] & 
& A_{PL}(E)\ar[rd] & 
\\ & 
(M\otimes _AN,d) \ar[rr]_{\hbox{\small$\sim$}} &
& A_{PL}(E\times_BE')
}$$
Form the pullback $(A,d)$--module $(J,D) = (M,d)\times _{(M\otimes _{A}N,d)}(Q,D)$ and the pullback cochain algebra $A_{PL}(E)\times _{A_{PL}(E\times_BE')}A_{PL}(Z)$. By the dual of the gluing lemma \cite[II.1.2]{Baues}, \cite[8.13]{GoerssJardine}, the horizontal quasi-isomorphisms in the above cube induce a quasi-isomorphism of $(A,d)$--modules  $$ (J,D) \we A_{PL}(E)\times _{A_{PL}(E\times_BE')}A_{PL}(Z).$$ By \cite[13.5]{RHT}, the canonical morphism $$A_{PL}(E*_BE') \to A_{PL}(E)\times _{A_{PL}(E\times_BE')}A_{PL}(Z)$$ is a quasi-isomorphism and we obtain the following commutative diagram of $(A,d)$--modules:
$$\disablesubscriptcorrection
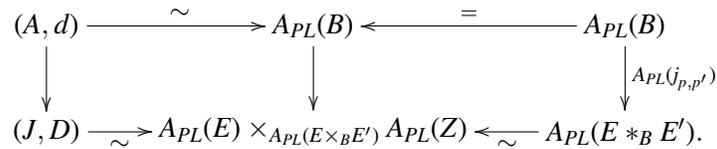
\xymatrix{
(A,d) \ar[r]^{\hbox{\small$\sim$}} \ar[d] 
& A_{PL}(B)\ar[d] 
& A_{PL}(B) \ar [d]^{A_{PL}(j_{p,p'} )} \ar[l]_{=} 
\\ 
(J,D) \ar[r]_-{\hbox{\small$\sim$}} 
& A_{PL}(E)\times _{A_{PL}(E\times_BE')}A_{PL}(Z)
& A_{PL}(E*_BE').\ar[l]^-{\hbox{\small$\sim$}} 
}$$

\subsection*{Step 3:\qua A quasi-isomorphism $(J,D) \we (M,d)\ast_{(A,d)}(N,d)$}  We have 
$$(J,D) = (M \oplus N \oplus s^{-1}M\otimes_{A}N,D).$$
The action on $M \oplus  N$ is the one of the direct sum, $a\cdot s^{-1}w = (-1)^{|a|}s^{-1}aw$, 
and the differential is given by $Dm = dm + s^{-1}\gamma(m)$, $Dn = dn + s^{-1}\nu(n)$, and $Ds^{-1}w = -s^{-1}dw$. Here, $\gamma \co  M \to M\otimes _AN$ is the inclusion. Let $j\co  (A,d) \to (J,D)$ be the canonical morphism. If we write $a_M$ for the elements of $J$ which lie in the copy of $A$ coming from $M$ and $a_N$ for the elements of $J$ which lie in the copy of $A$ coming from $N$ then $j$ is given by $j(a) = a_M - a_N$. 
Since the join $(M,d)\ast _{(A,d)}(N,d)$ is a semifree extension of $(A,d)$, by the lifting lemma \cite[II.1.11]{Baues}, in order to finish the proof it is enough to construct a quasi-isomorphism of $(A,d)$--modules $$f\co  (J,D) \to (M,d)\ast_{(A,d)}(N,d) = (A \oplus A\otimes s^{-1}X\otimes Y,d)$$ such that the composite of $f\circ j$ is the inclusion $(A,d) \to (M,d)\ast _{(A,d)}(N,d)$. We define the map $f$ by $f(a_M) = \disfrac{1}{2}a$, $f(a\otimes x) = 0$, $f(a_N) = -\disfrac{1}{2}a$, $f(a\otimes y) = 0$, $f(s^{-1}a) =0$, $f(s^{-1}a\otimes x) = -\disfrac{1}{2}(-1)^{|a|}ad_0x$, $f(s^{-1}a\otimes y) = \disfrac{1}{2}(-1)^{|a|}ad_0y$, and $f(s^{-1}a\otimes x \otimes y) = (-1)^{|a|}a\otimes s^{-1} x \otimes y$. It is straightforward to check that $f$ is $A$--linear and obvious that $f\circ j$ is the inclusion. Consider an element $x\in X$ and write $d_+x = \sum _ia_i\otimes x_i$. As we have shown at the beginning of the proof of \fullref{d^2=0}, $dd_0x = -\sum _i(-1)^{|a_i|}a_id_0x_i$. Using this identity and a  corresponding one for $y \in Y$, it is straightforward to check that $f$ commutes with the differentials. It remains to show that $f$ is a quasi-isomorphism. Consider the pushout $(A,d)$--module $(R,d) = (A \oplus A\otimes X \oplus A\otimes Y,d)$ of the inclusions $(A,d) \to (M,d)$ and $(A,d) \to (N,d)$ and form the acyclic differential vector space 
$$(R\oplus s^{-1}R,D) = ( A \oplus A\otimes X \oplus A\otimes Y \oplus s^{-1}A \oplus s^{-1}A\otimes X \oplus s^{-1}A\otimes Y,D)$$
where $Dr = dr + s^{-1}r$ and $Ds^{-1}r = -s^{-1}dr$. Define a map $g\co  R\oplus s^{-1}R \to J$
by $ga = \disfrac{1}{2}a_M + \disfrac{1}{2}a_N$, $g(a\otimes x) = a\otimes x$, $g(a\otimes y) = a\otimes y$,  $g(s^{-1}a) = s^{-1}a$, $$g(s^{-1}a\otimes x) = s^{-1}a\otimes x +\disfrac{1}{2}(-1)^{|a|}(ad_0x)_M - \disfrac{1}{2}(-1)^{|a|}(ad_0x)_N,$$ and
$$g(s^{-1}a\otimes y) = s^{-1}a\otimes y -\disfrac{1}{2}(-1)^{|a|}(ad_0y)_M + \disfrac{1}{2}(-1)^{|a|}(ad_0y)_N.$$
One easily checks that $f\circ g = 0$. Write $A_M$ to denote the copy of $A$ in $J$ coming from $M$. Then 
$$J = A_M \oplus \im g \oplus s^{-1}A\otimes X\otimes Y.$$
Therefore $g$ is an isomorphism onto $\ker f$. Using once more the identity $dd_0x = -\sum _i(-1)^{|a_i|}a_id_0x_i$, one checks that $g$ commutes with the differentials. Since  $(R\oplus  s^{-1}R,D)$ is acyclic, this implies that $f$ is a quasi-isomorphism. \hfill{$\Box$}

\section{The invariant Msecat}

\begin{defin}\rm
Let $p\co  E \to B$ be fibration. We define $\Msecat p$ to be the least integer $n$ such that there exists a commutative diagram of $A_{PL}(B)$--modules:
$$\disablesubscriptcorrection\xymatrix{
A_{PL}(B)\ar@{=}[r] \ar[d]_{A_{PL}(j^np)} \ar[dr]
& A_{PL}(B)
\\ 
A_{PL}(\ast^n_BE) 
& (P,d) \ar[l]^-{\hbox{\small$\sim$}} \ar[u] 
}$$
If no such $n$ exists we set $\Msecat p = \infty$.
\end{defin}

We first show that $\Msecat p$ is a lower bound of $\secat p$ which is closer to the sectional category than the classical lower bound $\nil \ker p^*$:

\begin{theor}\label{ineq}
For any fibration $p\co  E \to B$, $\nil \ker p^* \leq \Msecat p$. If $B$ is normal then $\Msecat p \leq \secat p$. 
\end{theor}

\dem
Suppose that $\Msecat p \leq n$. We show that $\nil \ker p^* \leq n$. Form the pullback: 
$$\disablesubscriptcorrection\xymatrix{
\ast^n_BE \times_BE \ar[r] \ar[d]_{\bar p_n} 
& E\ar [d]^{p} 
\\ 
\ast^n_BE \ar[r]_{j^np} 
& B
}$$
By \fullref{secatq}, $\secat \bar p_n \leq n$. Therefore $\nil \ker \bar p_n^* \leq n$. Since $\Msecat p \leq n$, the join map $j^np$ is injective in cohomology. Now consider elements $\alpha_0, \dots, \alpha_{n} \in \ker p^*$. Since $\nil \ker \bar p_n^* \leq n$, we have $(j^np)^*(\alpha_0\cup \cdots \cup \alpha_{n}) = 0$. Since $(j^np)^*$ is injective, $\alpha_0\cup \cdots \cup \alpha_{n} = 0$. This shows that $\nil \ker p^* \leq n$.

Suppose now that $B$ is normal and that $\secat p \leq n$. Then there exists a section $s\co  B \to \ast^n_BE$ of the join map $j^np$. We therefore have the following commutative diagram of $A_{PL}(B)$--modules:
$$\disablesubscriptcorrection\xymatrix{
A_{PL}(B)\ar@{=}[r] \ar[d]_{A_{PL}(j^np)} \ar[dr]
& A_{PL}(B)
\\ 
A_{PL}(\ast^n_BE) 
& A_{PL}(\ast^n_BE). \ar@{=}[l] \ar[u]_{A_{PL}(s)} 
}$$
It follows that $\Msecat p \leq n$. \findem

The number $\Msecat p$ can be calculated using the algebraic join construction of the previous sections. For the proof of this fact we need the following lemma:

\begin{lem} \label{retraction}
Let $p\co  E \to B$ be a fibration, $\alpha \co  (A,d) \we A_{PL}(B)$ be a commutative cochain algebra model, $i\co  (A,d) \to (Q,d)$ be a cofibration of $(A,d)$--modules, and $\phi \co  (Q,d) \to A_{PL}(\ast^n_BE)$ be a morphism of $(A,d)$--modules such that $\phi \circ i = A_{PL}(j^np) \circ \alpha$. If $\Msecat p \leq n$ then $i$ admits a retraction of $(A,d)$--modules.
\end{lem}

\dem
By definition, there is a commutative diagram of $A_{PL}(B)$--modules: $$\disablesubscriptcorrection\xymatrix{
A_{PL}(B)\ar@{=}[r] \ar[d]_{A_{PL}(j^np)} \ar[dr]^{j}
& A_{PL}(B)
\\ 
A_{PL}(\ast^n_BE) 
& (P,d) \ar[l]_-{\hbox{\small$\sim$}}^-{\psi} \ar[u]_{r} 
}$$
This is automatically a commutative diagram of $(A,d)$--modules. Form the following commutative diagram of $(A,d)$--modules:
$$\disablesubscriptcorrection\xymatrix{
(A,d) \ar[r]^{j\circ \alpha} \ar@{ >->}[d]_{i}
& (P,d) \ar [d]^{\psi}_{\hbox{\small$\sim$}} 
\\ 
(Q,d) \ar[r]_-{\phi}
& A_{PL}(\ast^n_BE)
}$$
By the lifting lemma \cite[II.1.11]{Baues}, there exists a morphism of $(A,d)$--modules $\lambda \colon  (Q,d)$ $ \to (P,d)$ such that $\lambda \circ i = j \circ \alpha$. We have obtained the following commutative diagram of $(A,d)$--modules:
$$\disablesubscriptcorrection\xymatrix{
(A,d) \ar@{=}[r] \ar@{ >->}[d]_{i} 
& (A,d)\ar [d]^{\hbox{\small$\sim$}}_{\alpha} 
\\ 
(Q,d) \ar[r]_{r\circ \lambda} 
& A_{PL}(B)
}$$
The lifting lemma \cite[II.1.11]{Baues} yields the required retraction of $i$. 
\findem

\begin{theor}\label{Mjoin}
Let $p\co  E \to B$ be a fibration between simply connected spaces of finite type, $\alpha \co  (A,d) \we A_{PL}(B)$ be a commutative cochain algebra model, and $(M,d) = (A\otimes (\Q \oplus X),d)$ be a semifree extension of $(A,d)$ such that there exists a quasi-isomorphism of $(A,d)$--modules $(M,d) \we A_{PL}(E)$ extending $A_{PL}(p)\circ \alpha$. Then $\Msecat p \leq n$ if and only if the inclusion $(A,d) \to \ast_{(A,d)}^n(M,d)$ admits a retraction of $(A,d)$--modules. 
\end{theor}

\dem
Suppose first that $\Msecat p \leq n$. By \fullref{joinmodel}, the morphism of\break $(A,d)$--modules $A_{PL}(j^np)\circ \alpha$ extends to a quasi-isomorphism of $(A,d)$--modules\break $\ast_{(A,d)}^n(M,d) \we  A_{PL}(\ast _B^nE)$. By  \fullref{retraction}, the inclusion $(A,d) \to \ast_{(A,d)}^n(M,d)$ admits a retraction of $(A,d)$--modules.

Suppose now that the inclusion $(A,d) \to \ast_{(A,d)}^n(M,d)$ admits a retraction $\rho$ of $(A,d)$--modules. Then the morphism of $A_{PL}(B)$--modules $A_{PL}(B)\otimes _{(A,d)}\rho$ is a retraction of the morphism of $A_{PL}(B)$--modules 
$$A_{PL}(B) = A_{PL}(B)\otimes_{(A,d)}(A,d) \to A_{PL}(B)\otimes_{(A,d)}\ast^n_{(A,d)}(M,d).$$
By \fullref{joinmodel}, the morphism of $(A,d)$--modules $A_{PL}(j^np)\circ \alpha$ extends to a quasi-isomor\-phism of $(A,d)$--modules $\psi\co  \ast_{(A,d)}^n(M,d) \we  A_{PL}(\ast _B^nE)$. Consider the following commutative diagram of $A_{PL}(B)$--modules: 
$$\disablesubscriptcorrection\xymatrix{
A_{PL}(B)\otimes_{A_{PL}(B)}A_{PL}(B) \ar[d]^{A_{PL}(B)\otimes_{A_{PL}(B)}A_{PL}(j^np)} 
&& A_{PL}(B)\otimes_{(A,d)}(A,d)\ar@{=}[ll]_{A_{PL}(B)\otimes _{\alpha}\alpha}  \ar[d] 
\\ 
A_{PL}(B)\otimes_{A_{PL}(B)}A_{PL}(\ast^n_BE)
&& A_{PL}(B)\otimes_{(A,d)}\ast^n_{(A,d)}(M,d) \ar[ll]^{A_{PL}(B)\otimes _{\alpha}\psi} 
}$$ 
Using \cite[6.10]{RHT} one sees that $A_{PL}(B)\otimes _{\alpha}\psi$ is a quasi-isomorphism. The left hand vertical morphism is precisely $A_{PL}(j^np)$. We have obtained the following commutative diagram of $A_{PL}(B)$--modules:
$$\disablesubscriptcorrection\xymatrix{
A_{PL}(B)\ar@{=}[rr] \ar[d]_{A_{PL}(j^np)} \ar[drr]
&& A_{PL}(B)
\\ 
A_{PL}(\ast^n_BE) 
&& A_{PL}(B)\otimes_{(A,d)}\ast^n_{(A,d)}(M,d) \ar[ll]_-{\hbox{\small$\sim$}}^-{A_{PL}(B)\otimes _{\alpha}\psi} \ar[u]_{A_{PL}(B)\otimes _{(A,d)}\rho} 
}$$
This shows that $\Msecat p \leq n$.
\findem

Note that we have not yet shown that $\Msecat $ is a homotopy invariant. This is contained in the following proposition.

\begin{prop}\label{Mbasics}
Consider a commutative diagram
$$\disablesubscriptcorrection\xymatrix{
E \ar[r]^{g} \ar[d]_{p} 
& E'\ar [d]^{p'} 
\\ 
B\ar[r]_{f} 
& B' 
}$$
in which $p$ and $p'$ are fibrations.

{\rm(a)}\qua If $f$ is a homotopy equivalence then $\Msecat p' \leq \Msecat p$.

{\rm(b)}\qua If $f$ and $g$ are homotopy equivalences then $\Msecat p' = \Msecat p$.

{\rm(c)}\qua If the diagram is a pullback and all spaces are simply connected and of finite type then $\Msecat p \leq \Msecat p'$.
\end{prop}

\dem
(a)\qua Suppose that $\Msecat p \leq n$. Choose a factorization $A_{PL}(j^np') = \psi \circ i$ where $i \co  A_{PL}(B') \to (Q,d)$ is a cofibration of $A_{PL}(B')$--modules and $\psi \co  (Q,d) \to A_{PL}(\ast^n_{B'}E')$ is a quasi-isomorphism of $A_{PL}(B')$--modules. Then $A_{PL}(\ast^n_fg) \circ \psi \circ i = A_{PL}(j^np)\circ A_{PL}(f)$. By \fullref{retraction}, $i$ admits a retraction of $A_{PL}(B')$--modules. This shows that $\Msecat p' \leq n$.

(b)\qua This is a formal consequence of (a). Indeed, by (a), $\Msecat p' \leq \Msecat p$. But if $f$ and $g$ are homotopy equivalences then the homotopy inverses can be used to construct a commutative square 
$$\disablesubscriptcorrection\xymatrix{
E'\ar[r]^{\simeq} \ar[d]_{p'} 
& E\ar[d]^{p} 
\\ 
B' \ar[r]_{\simeq} 
& B 
}$$
showing $\Msecat p \leq \Msecat p'$.

(c)\qua Applying the functor $A_{PL}$ to the given square we obtain the following commutative diagram of commutative cochain algebras:
$$\disablesubscriptcorrection\xymatrix{
A_{PL}(B')\ar[r]^{A_{PL}(f)} \ar[d]_{A_{PL}(p')} 
& A_{PL}(B)\ar [d]^{A_{PL}(p)} 
\\ 
A_{PL}(E')\ar[r]_{A_{PL}(g)} 
& A_{PL}(E)
}$$
Let $\alpha \co  (A',d) \we A_{PL}(B')$ be a Sullivan model. Choose factorizations  $A_{PL}(f)\circ \alpha = \psi \circ i$ and $A_{PL}(p')\circ \alpha = \phi \circ j$ such that 
$i \co  (A',d) \to (A,d)$ and $j \co  (A',d) \to (M',d)$ are inclusions of relative Sullivan algebras and $\psi \co  (A,d) \to A_{PL}(B)$ and $\phi \co  (M',d) \to A_{PL}(E')$ are quasi-isomorphisms. Then, by \cite[15(c)]{RHT}, the induced morphism of cochain algebras 
$$(M,d) = (A,d)\otimes _{(A',d)}(M',d) \to A_{PL}(E)$$ is a quasi-isomorphism. Note that $(M',d)$ is a semifree extension of $(A',d)$ and $(M,d)$ is a semifree extension of $(A,d)$.  
Suppose that $\Msecat p' \leq n$. By the preceding theorem, the inclusion $(A',d) \to \ast^n_{(A',d)}(M',d)$ admits a retraction $\rho$ of $(A',d)$--modules. As remarked in \fullref{iteratedjoin}(iii), $$\ast^n_{(A,d)}(M,d) = (A,d)\otimes _{(A',d)}\ast^n_{(A',d)}(M',d).$$ The morphism of $(A,d)$--modules 
$$A\otimes _{A'}\rho \co  (A,d)\otimes _{(A',d)}\ast^n_{(A',d)}(M',d)  \to (A,d)\otimes _{(A',d)}(A',d) = (A,d)$$
is a retraction of the inclusion $(A,d) \to \ast^n_{(A,d)}(M,d)$. By the preceding theorem, this implies that $\Msecat p \leq n$.
\findem

The next proposition shows that the invariant $\Msecat$ is a generalization of the well-known invariant $\Mcat$ of spaces. Let $B$ be a simply connected space of finite type with Sullivan model $(\Lambda V,d)$. By definition, $\Mcat B$ is the least integer $n$ such that for some (equivalently: any) Sullivan model $(\Lambda V \otimes \Lambda W,d) \we (\Lambda V / \Lambda ^{>n}V,d)$ of the projection $(\Lambda V,d) \to (\Lambda V / \Lambda ^{>n}V,d)$, the inclusion   $(\Lambda V,d) \to (\Lambda V \otimes \Lambda W,d)$ admits a retraction of $(\Lambda V,d)$--modules.
If no such $n$ exists, $\Mcat B = \infty$.

\begin{prop}\label{Mcat}
Let $B$ be a simply connected pointed space of finite type. Consider the evaluation fibration $ev_1\co  PB \to B$, $\omega \mapsto \omega(1)$. Then $\Msecat ev_1 = \Mcat B$.
\end{prop}

\dem Let $\alpha \co (\Lambda V,d) \we A_{PL}(B)$ be a Sullivan model
of $B$. Denote the projection $(\Lambda V,d) \to (\Lambda V/\Lambda
^{>n}V,d)$ by $q_n$ and choose a factorization $q_n = \phi \circ i$
where $i \co (\Lambda V,d) \to (\Lambda V\otimes \Lambda W,d)$ is the
inclusion of a relative Sullivan algebra and $\phi \co (\Lambda
V\otimes \Lambda W,d)\to (\Lambda V/\Lambda ^{>n}V,d)$ is a
quasi-isomorphism. Choose a factorization $A_{PL}(j^nev_1)\circ \alpha
= \psi \circ j$ where $j \co (\Lambda V,d) \to (\Lambda V\otimes
\Lambda X,d)$ is the inclusion of a relative Sullivan algebra and
$\psi \co (\Lambda V\otimes \Lambda X,d)\to A_{PL}(\ast^n_BPB)$ is a
quasi-isomorphism. It follows from F\'elix and Halperin \cite{FH} that
there exist morphisms of commutative cochain algebras $\sigma \co
(\Lambda V\otimes \Lambda X,d) \to (\Lambda V\otimes \Lambda W,d)$ and
$\rho\co (\Lambda V\otimes \Lambda W,d) \to (\Lambda V\otimes \Lambda
X,d)$ such that $\sigma \circ j = i$ and $\rho \circ i = j$. This
implies that $\Mcat B \leq n$ if and only if $j$ admits a retraction
of $(\Lambda V,d)$--modules. Let $(M,d)$ be a semifree extension of
$(\Lambda V,d)$ such that there exists a quasi-isomorphism of
$(\Lambda V,d)$--modules $(M,d) \we A_{PL}(PB)$ which extends
$A_{PL}(ev_1) \circ \alpha$. By \fullref{joinmodel}, the morphism
of $(\Lambda V,d)$--modules $A_{PL}(j^nev_1)\circ \alpha$ extends to a
quasi-isomorphism of $(\Lambda V,d)$--modules $\ast_{(\Lambda
V,d)}^n(M,d) \we A_{PL}(\ast _B^nPB)$. Use the lifting lemma
\cite[II.1.11]{Baues} to construct quasi-isomorphisms of $(\Lambda
V,d)$--modules $\beta \co \ast^n_{(\Lambda V,d)}(M,d) \we (\Lambda
V\otimes \Lambda X,d)$ and $\gamma \co (\Lambda V\otimes \Lambda X,d)
\we \ast^n_{(\Lambda V,d)}(M,d)$ such that $\beta$ extends $j$ and
$\gamma \circ j$ is the inclusion of $(\Lambda V,d)$. We obtain that
$\Mcat B \leq n$ if and only if the inclusion $(\Lambda V,d) \to
\ast^n_{(\Lambda V,d)}(M,d)$ admits a retraction of $(\Lambda
V,d)$--modules. By \fullref{Mjoin}, this is the case if and only if
$\Msecat ev_1 \leq n$. \findem

Recall that an upper bound for the sectional category of a surjective fibration is given by the Lusternik--Schnirelmann category of the base space. The following is the analogous result for $\Msecat$ and $\Mcat$.

\begin{prop}\label{MsecatleqMcat}
Let $B$ be a simply connected space of finite type. For any surjective fibration $p \co  E \to B$, $\Msecat p \leq \Mcat B$.
\end{prop}

\dem
Recall that for us simply connected spaces are non-empty. Fix any base point in $B$ and consider the evaluation fibration $ev_1 \co  PB \to B$. Since $p$ is surjective, $E\not= \emptyset$. Since $PB$ is contractible, there exists a continuous map $\lambda \co  PB \to E$ such that $p\circ \lambda = ev_1$. By \fullref{Mbasics}(a) and \fullref{Mcat}, $\Msecat p \leq \Msecat ev_1 = \Mcat B$.
\findem

\begin{rem}\rm
K Hess \cite{Hess} has shown that the invariant $\Mcat$ coincides for simply connected spaces with rational Lusternik--Schnirelmann category. This result does not generalize to sectional category. Indeed, D Stanley \cite{Stanley} has constructed a fibration $p$ with fiber $S^2$ whose rational sectional category is $1$. By Vandembroucq \cite{Transfert}, any such fibration satisfies $\Msecat p = 0$.
\end{rem}\rm

\section{Topological complexity}
In \cite{Farber1} and  \cite{Farber2}, M Farber defined the \textit{topological complexity} of a space $X$, $\TC(X)$, to be the sectional category of the evaluation fibration $ev_{0,1} \co  X^I \to X\times X, \omega \mapsto (\omega (0), \omega (1))$. This invariant has proved to be very useful in the study of the motion planning problem in robotics. Note that Farber's definition of $\TC$ differs by $1$ from the one given here. In this section we study the invariant 
$$\MTC(X) = \Msecat (ev_{0,1} \co  X^I \to X\times X).$$
In order to simplify the presentation we restrict our attention to simply connected spaces of finite type having the homotopy type of CW complexes. 

The evaluation fibration $ev_{0,1} \co  X^I \to X\times X$ is the mapping path fibration associated to the diagonal map $\Delta \co  X \to X\times X$. We may therefore identify the map $ev^*_{0,1} \co  H^*(X \times X) \to H^*(X^I)$ with the cup product $\cup \co  H^*(X)\otimes H^*(X) \to H^*(X)$. 

\begin{prop}\label{MTCineq}
We have $\nil \ker \cup \leq \MTC(X) \leq \TC(X)$ and $\Mcat X \leq  \MTC(X)$ $\leq 2\Mcat X$. 
\end{prop}

\dem
The first inequalities follow from \fullref{ineq}. By \fullref{MsecatleqMcat} and \cite[30.2]{RHT}, $\MTC(X) \leq \Mcat(X\times X) = 2\Mcat X$. For the remaining inequality consider the map $f\co  X \to X\times X$, $x \mapsto (*,x)$ where $* \in X$ is any base point and form the following pullback diagram:
$$\disablesubscriptcorrection\xymatrix{
PX \ar[r] \ar[d]_{ev_1} 
& X^I\ar [d]^{ev_{0,1}} 
\\ 
X \ar[r]_-{f} 
& X\times X 
}$$
By \fullref{Mcat} and \fullref{Mbasics}(c), $\Mcat X = \Msecat ev_1 \leq \Msecat ev_{0,1} = \MTC(X)$. 
\findem

Consider a space $X$ with Sullivan model $(\Lambda V,d)$. A Sullivan model of the product space $X\times X$ is then given by
$(\Lambda(V\oplus V'),d) = (\Lambda V,d)\otimes (\Lambda V',d)$ where $(\Lambda V',d)$ is second copy of $(\Lambda V,d)$.
As is shown in \cite[pages 206--207]{RHT}, a model of the evaluation fibration (and the diagonal map) is given by the inclusion
$$(\Lambda (V\oplus V'),d) \to (\Lambda (V \oplus V')\otimes  \Lambda \overline{V},d)$$
where 
$$d(\bar v) = v' - v - \sum \limits _{i=1}^{\infty} \disfrac{(\zeta d)^i}{i!}(v).$$
Here, $\zeta$ is the derivation of degree $-1$ defined by $\zeta (v) = \zeta(v') = \bar v$ and $\zeta ( \bar v) = 0$. Using this explicit semifree extension of $(\Lambda (V\oplus V'),d)$ and the formula for the differential of the iterated join given in \fullref{iteratedjoin}(ii), one can calculate the invariant $\MTC(X)$ from $(\Lambda V,d)$. We remark that $d_0\bar v = v' -v$ and that $d_0x = 0$ for $x \in \Lambda ^{>1}\overline{V}$. Note also that if $(\Lambda V,d)$ is the minimal Sullivan model of $X$ then $(\Lambda (V \oplus V')\otimes  \Lambda \overline{V},d)$ is a minimal semifree $(\Lambda (V \oplus V'),d)$--module. 

The following proposition provides an upper bound for $\MTC$.

\begin{prop}\label{nillemma}
Let $(A,d)$ be a commutative cochain algebra model of $X$ with multiplication $\mu$. Then $\MTC(X) \leq \nil \ker \mu$.
\end{prop}

\dem
Suppose that $\nil \ker \mu \leq n$. We show that the inclusion
$$i\co  (\Lambda (V\oplus V'),d) \to \ast^n_{(\Lambda (V\oplus V'),d)}(\Lambda (V \oplus V')\otimes  \Lambda \overline{V},d)$$
admits a retraction of $(\Lambda (V\oplus V'),d)$--modules. Choose a quasi-isomorphism of commutative cochain algebras $\alpha \co  (\Lambda V,d) \we (A,d)$. Consider the tensor product algebra $(A,d)\otimes (A,d) = (A\otimes A,d)$ and the $(A\otimes A,d)$--semifree extension of $(A\otimes A,d)$ defined by
$$(M,d) = (A\otimes A,d)\otimes _{(\Lambda (V\oplus V'),d)}(\Lambda (V \oplus V')\otimes  \Lambda \overline{V},d) = (A\otimes A\otimes \Lambda \overline{V},d).$$
We have 
$$\ast^n_{(A\otimes A,d)}(M,d) = (A\otimes A \oplus A\otimes A \otimes s^{-n}(\Lambda^+\overline{ V})^{\otimes n+1},d).
$$
Consider an element $s^{-n}x_0\otimes \cdots \otimes x_n \in s^{-n}(\Lambda^+\overline{ V})^{\otimes n+1}$. If one of the $x_i$ lies in $\Lambda ^{>1}\overline{V}$ then $d(s^{-n}x_0\otimes \cdots \otimes x_n)$ has no term in $A\otimes A$. 
Since $\nil \ker \mu \leq n$, this also holds if all $x_i \in \overline{ V}$. We can thus define an $(A\otimes A,d)$--module retraction $r$ of the inclusion  $j\co  (A\otimes A,d) \to \ast^n_{(A\otimes A,d)}(M,d)$ by sending  $A \otimes A \otimes s^{-n}(\Lambda^+\overline{ V})^{\otimes n+1}$ to $0$. By \fullref{iteratedjoin}(iii), the map $j$ is obtained by applying the functor $(A\otimes A,d)\otimes _{(\Lambda (V\oplus V'),d)}-$ to the inclusion $i$. Consider the following commutative diagram of $(\Lambda (V\oplus V'),d)$--modules:
$$\disablesubscriptcorrection\xymatrix{
(\Lambda (V\oplus V'),d) \ar@{=}[rr] \ar@{ >->}[d]_{i} 
& &(\Lambda (V\oplus V'),d) \ar [d]^{\hbox{\small$\sim$}}_{\alpha \otimes \alpha} 
\\ 
\ast^n_{(\Lambda (V\oplus V'),d)}(\Lambda (V \oplus V')\otimes  \Lambda \overline{V},d) \ar[rr]_-{r\circ ((\alpha \otimes \alpha) \otimes_{id}id)} 
& &(A\otimes A,d)
}$$
The lifting lemma \cite[II.1.11]{Baues} yields the required retraction of $i$. 
\findem

Note that the number $\nil \ker \mu$ is in general not the same for different commutative cochain algebra models of $X$. Given a commutative graded algebra $A$ with multiplication $\mu$, the number $\nil \ker \mu$ can be determined using the following lemma:

\begin{lem}  \label{nillemma2}
Let $\Lambda W$ be a commutative graded algebra and $I \subset  \Lambda W$ be an ideal such that $A = \Lambda W/I$.
Let $\Lambda W'$ be a second copy of $\Lambda W$ and  $J  \subset  \Lambda (W\oplus W')= \Lambda W \otimes \Lambda W'$ be the
ideal $I \otimes \Lambda W' + \Lambda W\otimes I'$ where $I'$ is the ideal of $\Lambda W'$ corresponding to $I$. Let finally $\mathcal B$
be a basis of the graded vector space $W$. Then $\nil \ker \mu \leq n$ if and only if, for all $w_0, \dots , w_n \in \mathcal B$,
$(w_0' -w_0)\cdots (w_n' -w_n) \equiv 0 \mod J$.
\end{lem}

\dem
Denote the multiplication $\Lambda (W\oplus W') = \Lambda W\otimes \Lambda W' \to \Lambda W$ by $m$. We have
$A\otimes A = \Lambda (W\oplus W')/J$ and $\ker \mu = \ker m /(J \cap \ker m)$. It suffices to show that $\ker m$ is the ideal of
$\Lambda (W\oplus W')$ generated by the elements $w'-w$, $w \in \mathcal B$. Denote this ideal by $K$. Obviously, $K\subset \ker m$.
Note also that $K\cap \Lambda W = 0$. In order to show the equality $K = \ker m$ we show that $\Lambda W \oplus K =
\Lambda(W\oplus W')$. For this it is enough to show that for each $n\geq 1$, $\Lambda W\otimes \Lambda ^{n}W' \subset  \Lambda
W \oplus K$. We proceed by induction. Consider $w\in \mathcal B$ and $\xi \in \Lambda W$. We have $\xi w' = \xi w + \xi (w'-w)
\in \Lambda W \oplus K$. Suppose the assertion holds for some $n \geq 1$. Consider $\zeta \in  \Lambda W\otimes \Lambda ^{n}W'$
and $w\in \mathcal B$. By the inductive hypothesis, $\zeta \in \Lambda W \oplus K$. Write $\zeta = \theta + k$ with $\theta \in \Lambda W$
and $k \in K$. By the inductive hypothesis,
$\theta w' \in \Lambda W \oplus K$. Since $K$ is an ideal, $kw' \in K$. It follows that $\zeta w' \in \Lambda W \oplus K$.  This closes the
induction and the result follows.
\findem

Recall that a space $X$ is called \textit{formal} if $H^*(X)$ is a commutative cochain algebra model of $X$.  \fullref{nillemma} immediately implies the following:

\begin{prop}
If $X$ is formal then $\MTC(X) = \nil \ker \cup$.
\end{prop}

\begin{exam}\rm
The simplest example of a non-formal space is the space $$X = S^3_a \vee S^3_b \cup e^8 \cup e^8$$ where the
$8$--cells are attached by means of the iterated Whitehead products\break  $[S^3_a, [S^3_a,S^3_b]]$ and $[S^3_b, [S^3_a,S^3_b]]$. We show that this space satisfies $\MTC(X)= 3$ and $\nil \ker \cup = 2$. For degree reasons, the space $X$ has the same cohomology algebra as the wedge of spheres $S^3\vee S^3 \vee S^8 \vee S^8$. Therefore $X$ satisfies $\nil \ker \cup = 2$. Indeed, since $\TC \leq 2\cat$, the topological complexity of a wedge of spheres is $\leq 2$. On the other hand, any space with at least two cohomology generators satisfies $\nil \ker \cup \geq 2$: if $\xi$ and $\theta$ are two cohomology generators then $(\xi \otimes 1 - 1\otimes \xi)(\theta \otimes 1 -1 \otimes \theta)$ is a non-trivial product in $\ker \cup$. The minimal Sullivan model of $X$ is the algebra $(\Lambda(V),d)$ where the graded $\Q$--vector space $V$ is generated by cocycles $a$ and $b$ of degree $3$, an element $u$ of degree $5$ with $du = ab$, and elements of degree $> 8$. Consider the $d$--stable ideal $I = (\Lambda V)^{\geq 9}$ and form the quotient algebra $(A,d) = (\Lambda V/I,d)$. Since $(I,d)$ is acyclic, the projection $(\Lambda V,d) \to (\Lambda V/I,d)$ is a quasi-isomorphism. Consider the ideal $J = I \otimes \Lambda V' + \Lambda V \otimes I' \subset \Lambda (V \oplus V')$ as in \fullref{nillemma2}.
We have $(a'-a)(b'-b)(u'-u)\not \equiv 0\mod J$. Since $a,b, u$ are of odd degree, any longer non-zero 
product of the form $(v_0' -v_0)\cdots (v_n' -v_n)$ must contain at least one factor $v_i' - v_i$ with $|v_i| \geq 9$. 
For $n \geq 3$ any such product is therefore an element of $J$. By \fullref{nillemma} and \fullref{nillemma2}, this implies that $\MTC(X) \leq 3$. 
We show that $\MTC(X) > 2$. The differential of the generators $\bar a$, $\bar b$, and $\bar u$ of the model $(\Lambda (V \oplus V')\otimes\Lambda \overline{V},d)$ of $X^I$ is given by $d\bar a = a' - a$, $d\bar b = b' -b$, and $d\bar u = u' - u + \alpha\otimes \bar a + \beta \otimes \bar b$ where $\alpha, \beta \in \Lambda (V\oplus V')$ are some elements of degree $3$. A straightforward calculation shows that $\alpha = -\frac{1}{2}(b+b')$ and $\beta = \frac{1}{2}(a+a')$. This information is, however, not needed for the calculations. It suffices to show that the inclusion
\begin{eqnarray*}
(\Lambda (V\oplus V'),d) &\to& \ast^2_{(\Lambda (V\oplus V'),d)}(\Lambda (V\oplus V') \otimes \Lambda \overline{V},d)\\
&=& (\Lambda (V \oplus V')\oplus \Lambda (V \oplus V')\otimes s^{-2}\Lambda ^+\overline{V}\otimes \Lambda ^+\overline{V}\otimes \Lambda ^+\overline{V},d)
\end{eqnarray*}
is not injective in cohomology. The element $z = (a'-a)(b'-b)(u'-u)$ is a cocycle of degree $11$ in $\Lambda(V\oplus V')$ which is not a coboundary. In the 2--fold join, however, we have
$$d(s^{-2}(\bar a\otimes \bar b\otimes \bar u + \bar b\otimes \bar u\otimes \bar a + \bar u\otimes \bar a\otimes \bar b - \bar a\otimes \bar u\otimes \bar b - \bar b\otimes \bar a\otimes \bar u - \bar u\otimes \bar b\otimes \bar a)) = -6z$$ so that $[z] = 0 \in H^{11}(\ast^2_{(\Lambda (V\oplus V'),d)}(\Lambda (V\oplus V') \otimes \Lambda \overline{V},d))$.
\end{exam}

Our last result is the fact that the difference $\MTC(X) - \nil \ker \cup$ can be arbitrarily large:

\begin{prop}
{\rm(i)}\qua For any $n \in \N$ there exists a finite CW--complex $X$ such that $\MTC(X) - \nil \ker \cup \geq n$.

{\rm(ii)}\qua There exists a space $X$ such that $\MTC(X) = \infty$ and $\nil \ker \cup < \infty$.
\end{prop}

\dem
(i)\qua Let $Z$ be a simply connected finite CW--complex having the same cohomology algebra as a wedge of spheres $Y$ and satisfying $\Mcat Z = 3$. Such a space has for instance been constructed by Kahl and Vandembroucq \cite{gaps}. Let $X$ be the $n$--fold product of the space $Z$, $X = Z^n$. Then $X$ is a finite CW--complex which has the same cohomology algebra as $Y^n$ and satisfies $\nil \ker \cup \leq \TC(Y^n) \leq 2\cat (Y^n) \leq 2n \cat Y = 2n$. On the other hand, by \fullref{MTCineq} and \cite[30.2]{RHT}, $\MTC(X) \geq \Mcat (X) = n \Mcat (Z) = 3n$.

(ii)\qua It suffices to take a space $X$ such that $\Mcat X = \infty$ and $\nil H^+(X) < \infty$. Such a space has been constructed by F\'elix, Halperin and Thomas \cite{FHT83}.
\findem

\bibliographystyle{gtart}
\bibliography{link}

\begin{thebibliography}{}
\providecommand\bibmarginpar{\leavevmode\marginpar}
\def\urlstyle#1{{\tt #1}}

\bibitem{Baues}
\textbf{H\,J Baues}, \emph{Algebraic homotopy}, Cambridge Studies in Advanced
  Mathematics 15, Cambridge University Press, Cambridge (1989) \xox{MR}{985099}

\bibitem{CLOT}
\textbf{O Cornea}, \textbf{G Lupton}, \textbf{J Oprea}, \textbf{D Tanr{\'e}},
  \emph{Lusternik--{S}chnirelmann category}, Mathematical Surveys and
  Monographs 103, American Mathematical Society, Providence, RI (2003)
  \xox{MR}{1990857}

\bibitem{Farber1}
\textbf{M Farber}, \href{http://dx.doi.org/10.1007/s00454-002-0760-9}
  {\emph{Topological complexity of motion planning}}, Discrete Comput. Geom. 29
  (2003) 211--221 \xox{MR}{1957228}

\bibitem{Farber2}
\textbf{M Farber}, \href{http://dx.doi.org/10.1016/j.topol.2003.07.011}
  {\emph{Instabilities of robot motion}}, Topology Appl. 140 (2004) 245--266
  \xox{MR}{2074919}

\bibitem{Agnese}
\textbf{A Fass\`o~Velenik}, \emph{Relative homotopy invariants of the type of
  the Lusternik-Schnirelmann category}, PhD thesis, FU Berlin (2002)

\bibitem{FH}
\textbf{Y F{\'e}lix}, \textbf{S Halperin},
  \href{http://links.jstor.org/sici?sici=0002-9947(198209)273:1%3C1:RLCAIA%3E2%
.0.CO%3B2-9} {\emph{Rational {LS} category and its applications}}, Trans. Amer.
  Math. Soc. 273 (1982) 1--38 \xox{MR}{664027}

\bibitem{FHT83}
\textbf{Y F{\'e}lix}, \textbf{S Halperin}, \textbf{J-C Thomas},
  \href{http://www.numdam.org/item?id=BSMF_1983__111__89_0} {\emph{L.{S}.
  cat\'egorie et suite spectrale de {M}ilnor--{M}oore (une nuit dans le
  train)}}, Bull. Soc. Math. France 111 (1983) 89--96 \xox{MR}{710377}

\bibitem{FHTDGA}
\textbf{Y F{\'e}lix}, \textbf{S Halperin}, \textbf{J-C Thomas},
  \emph{Differential graded algebras in topology}, from: ``Handbook of
  algebraic topology'', North-Holland, Amsterdam (1995)  829--865
  \xox{MR}{1361901}

\bibitem{RHT}
\textbf{Y F{\'e}lix}, \textbf{S Halperin}, \textbf{J-C Thomas}, \emph{Rational
  homotopy theory}, Graduate Texts in Mathematics 205, Springer, New York
  (2001) \xox{MR}{1802847}

\bibitem{GoerssJardine}
\textbf{P\,G Goerss}, \textbf{J\,F Jardine}, \emph{Simplicial homotopy theory},
  Progress in Mathematics 174, Birkh\"auser Verlag, Basel (1999)
  \xox{MR}{1711612}

\bibitem{Hess}
\textbf{K\,P Hess}, \href{http://dx.doi.org/10.1016/0040-9383(91)90006--P}
  {\emph{A proof of {G}anea's conjecture for rational spaces}}, Topology 30
  (1991) 205--214 \xox{MR}{1098914}

\bibitem{James1}
\textbf{I\,M James}, \href{http://dx.doi.org/10.1016/0040-9383(78)90002-2}
  {\emph{On category, in the sense of {L}usternik--{S}chnirelmann}}, Topology
  17 (1978) 331--348 \xox{MR}{516214}

\bibitem{James2}
\textbf{I\,M James}, \emph{Lusternik--{S}chnirelmann category}, from:
  ``Handbook of algebraic topology'', North-Holland, Amsterdam (1995)
  1293--1310 \xox{MR}{1361912}

\bibitem{gaps}
\textbf{T Kahl}, \textbf{L Vandembroucq},
  \href{http://projecteuclid.org/getRecord?id=euclid.bbms/1102715104}
  {\emph{Gaps in the {M}ilnor--{M}oore spectral sequence}}, Bull. Belg. Math.
  Soc. Simon Stevin 9 (2002) 265--277 \xox{MR}{2017081}

\bibitem{Svarc}
\textbf{A Schwarz}, \emph{The genus of a fiber space}, Amer. Math. Soc. Transl.
  55 (1966) 49--140

\bibitem{Smale}
\textbf{S Smale}, \emph{On the topology of algorithms I}, J. Complexity 3
  (1987) 81--89 \xox{MR}{907191}

\bibitem{Stanley}
\textbf{D Stanley},
  \href{http://www.ams.org/jourcgi/jour-getitem?pii=S000299390005468X}
  {\emph{The sectional category of spherical fibrations}}, Proc. Amer. Math.
  Soc. 128 (2000) 3137--3143 \xox{MR}{1691006}

\bibitem{Transfert}
\textbf{L Vandembroucq},
  \href{http://projecteuclid.org/getRecord?id=euclid.bbms/1103055585}
  {\emph{Sur le transfert de {B}ecker et {G}ottlieb}}, Bull. Belg. Math. Soc.
  Simon Stevin 6 (1999) 605--614 \xox{MR}{1732893}

\end{thebibliography}

\end{document}